\documentstyle{article}

\newcommand{\ev}{{\rm ev}}
\newcommand{\R}{{\bf R}}

\newcommand{\Q}{{\bf Q}}

\newcommand{\Z}{{\bf Z}}
\newcommand{\C}{{\bf C}}

\newcommand{\Hh}{{\cal H}}

\newcommand{\Ll}{{\cal L}}

\newcommand{\oMm}{{\overline{\cal M}}}

\newcommand{\id}{{\rm id}}
\newcommand{\Tilde}{\widetilde}
\newcommand{\p}{{\partial}}
\newcommand{\al}{{\alpha}}

\newcommand{\be}{{\beta}}

\newcommand{\Om}{{\Omega}}
\newcommand{\om}{{\omega}}

\newcommand{\Ga}{{\Gamma}}
\newcommand{\ka}{{\kappa}}
\newcommand{\la}{{\lambda}}
\newcommand{\La}{{\Lambda}}
\newcommand{\si}{{\sigma}}

\newcommand{\Symp}{{\rm Symp}}
\newcommand{\Diff}{{\rm Diff}}
\newcommand{\Ham}{{\rm Ham}}
\newcommand{\Hom}{{\rm Hom}}

\newcommand{\PD}{{\rm PD}}

\newcommand{\Si}{{\Sigma}}

\newcommand{\MS}{{\medskip}}
\newcommand{\SS}{{\smallskip}}

\newcommand{\NI}{{\noindent}}
\newcommand{\proof}[1]{\noindent{\bf Proof#1:\  }}
\newcommand{\QED}{\hfill$\Box$\medskip}

\begin{document}

\title{Topological rigidity of Hamiltonian
loops and quantum homology}
\author{Fran\c{c}ois Lalonde\thanks{Partially supported by NSERC grant 
OGP 0092913
and FCAR grant ER-1199.} \\ Universit\'e du Qu\'ebec \`a Montr\'eal
\\ (flalonde@math.uqam.ca) \and Dusa McDuff\thanks{Partially
supported by NSF grant DMS 9704825.} \\ State University of New York
at Stony Brook \\ (dusa@math.sunysb.edu)
 \and Leonid
Polterovich\thanks{All three authors supported by a NSERC
Collaborative Project Grant CPG0163730.} \\
Tel-Aviv University \\ (polterov@math.tau.ac.il)}

\date{October 1, 1997}
\maketitle

\NI
{\large \bf\S1. Main result}

\medskip

Let $G$ be a connected subgroup of the group  
$\Diff(M)$ of diffeomorphisms
of a manifold $M$. It is well known that every element
$\phi \in \pi_1(G,\id)$ defines an endomorphism
$\partial _{\phi} : H_*(M,\Q) \to H_{*+1}(M,\Q)$ as follows.
Choose a loop $\{\phi_t\},\; t\in S^1,$ of diffeomorphisms from $G$
representing $\phi$ and a cycle $C$ in $M$. Then the homology
class $\partial_{\phi} ([C])$ is represented by the cycle 
$S^1\times C\to M$ which is
spanned by 
$C$ under the loop $\{\phi_t\}$.

Suppose now that $(M,\omega)$ is a closed symplectic manifold,
and take $G$ to be  its group 
$\Ham(M,\omega)$ of Hamiltonian
diffeomorphisms. 
In this paper, we discuss the following statement. As
we will see in section 2 it has a number of applications
to the geometry and topology of the group of symplectomorphisms.

 \proclaim Theorem 1.A.
Let  $\phi$ be a loop in
the group $\Ham(M,\om)$ of Hamiltonian diffeomorphisms. Then
$\partial_{\phi}$ vanishes identically for all $\phi \in
\pi_1(\Ham(M,\om),\id)$.


Below we give the proof of this statement when $M$
is $4$\/-dimensional as well
as for some higher dimensional symplectic manifolds
-- the so-called spherically monotone manifolds and other manifolds 
 where $J$-holomorphic curves are well-behaved. See \S6 for the precise
assumptions on $(M, \om)$.
We believe that all theorems stated in this paper do hold in
the general case too, which forms the subject of our forthcoming paper
 [LMP2]; it
requires the use of the 
Gromov--Witten invariants on arbitrary manifolds
as developed by Fukaya--Ono~[FO], Li--Tian~[LiT], 
Liu--Tian~[LT], Ruan [R], and Hofer-Salamon [HS], and involves a deeper 
and more technical analysis. There are still a number of
details to be checked.

However, all the basic geometric ideas are 
already present in the particular case treated here.
An early version of these ideas was described in the survey article
by McDuff~[M2].

These ideas were inspired by a recent result of Seidel [Se] 
who discovered  a canonical action of a certain extension
of the group $\pi_1(\Ham(M,\om))$ on the quantum homology ring of $M$
that   arises from the natural action of the element $\phi$
on the loop space of $M$.  Seidel defines this action under
the additional assumptions mentioned above, and we will show that 
in this case 1.A  can be deduced from Seidel's result by simple
geometric arguments.

Notice that the particular case of 1.A  stating that the map
$$
\p_\phi: H_0(M,\Q)
\to H_1(M,\Q)
$$
 vanishes simply means  that the orbits of a periodic 
Hamiltonian flow \footnote{In this note, when no specific mention is made,
``Hamiltonian flow'' is understood
in the general sense of a time-dependent Hamiltonian.}  
are homologous to zero.  This  is a classical result and is very
easy to prove (see [BP], II-1.3).  Also when $\phi$ is a 
Hamiltonian circle action the
statement of 1.A immediately follows from a result of Kirwan (see 3.C below.)

  To our knowledge, the results of the present paper constitute
the first application of Quantum homology to Hamiltonian mechanics.

 The paper is organized as follows. 
In \S 2 we explain what  ``topological rigidity"
is, and  we derive from 1.A a new case
of the Flux Conjecture.
In \S 3 we reformulate 1.A in the more geometric
language of symplectic fibrations over the 2-sphere.
In \S 4 we describe Seidel's action and give the
proof of 1.A. in the particular cases explained in \S6. 
\S 5 contains a refinement of Seidel's theory.  In particular 
we construct a representation of the 
group $\pi_1(\Ham(M,\om))$ into the group
of automorphisms of an extension of the usual quantum 
cohomology ring of $M$. 
Finally, in the last section, \S6, we discuss the conditions on $M$
under which the results of the present paper hold, and what is
therefore left to prove in [LMP2].


\bigskip

\NI
{ \large\bf \S2. Rigidity of Hamiltonian loops}

\medskip

Let $\Symp_0(M,\om)$ be
the connected component of the identity
in the group of all symplectomorphisms
of $(M,\om)$.  We will say that the class $\phi\in \pi_1(\Diff(M), \id)$
has an $\om$-symplectic representative if it may be represented by a loop
$\{\phi_t\},\; t \in S^1,$  
in  $\Symp_0(M,\om)$
based at the identity.  The loop $\{\phi_t\}$ is $\om$-Hamiltonian if it
is the flow (with respect to $\om$) of a time-dependent Hamiltonian function
$H_t$.  We are interested in the question of which classes in
$\pi_1(\Symp_0(M,\om),\id)$ have Hamiltonian representatives, or equivalently
of when $\{\phi_t\}$ is
homotopic (through $\om$-symplectic loops) to a Hamiltonian loop.

\proclaim Theorem 2.A.  Suppose that $\om_1$ and $\om_2$ are two symplectic
forms on $M$ and that $\phi\in \pi_1(\Diff(M), \id)$ contains $\om_i$-symplectic
representatives $\{\phi_t^i\}$, for $i = 1,2$.  Then $\phi$ contains a
$\om_1$-Hamiltonian representative if and only if it contains a 
$\om_2$-Hamiltonian representative.

In other words, once we know that a 
loop has a symplectic representative,
the question of whether or not this representative can be chosen 
to be Hamiltonian
is independent of the choice of the symplectic form, 
and in particular of its  cohomology
class. This forms the content of the phenomenon of 
{\it topological rigidity of 
Hamiltonian loops} which is announced in the 
title of the present note. \footnote{
Actually, the proof of  Theorem 2.A shows that the following
stronger statement holds:  if the class
$\phi$ contains a $\om_1$-Hamiltonian representative, then any 
$\om_2$-symplectic loop is homotopic through $\om_2$-symplectic loops to a
$\om_2$-Hamiltonian loop. This result implies the following slighly
different version of 2.A: the loop $\{\phi_t^1\}$ is
homotopic (through $\om_1$-symplectic loops)
to a Hamiltonian loop if and only if the same is true for $\{\phi_t^2\}$.}

One situation in which this applies is when $\om'$
is a $C^{\infty}$-small perturbation of $\om$ in the space of closed $2$-forms.
Using Moser's argument one can easily show that  any given 
$\om$-symplectic loop $\{\phi_t\}$ can be perturbed to an
$\om'$-symplectic loop $\{\phi_t'\}$ provided that $\om'$ is
sufficiently $C^\infty$-close to $\om$.

\proclaim Corollary 2.B. Assume in the above situation that 
the loop $\{\phi_t\}$ is homotopic  in $\Symp_0(M,\om)$
to a Hamiltonian loop.
Then the loop $\{\phi_t'\}$ is homotopic  in $\Symp_0(M,\om')$
to a Hamiltonian loop.

In other words, the property of a loop of symplectomorphisms
to be Hamiltonian up to homotopy is stable with respect to
(small)  deformations of the symplectic structure.

The above theorem is an almost immediate consequence of 1.A
because of the characterization of Hamiltonian loops via the flux
homomorphism. Recall that the {\it flux homomorphism} 
$$
F_\om: \pi_1(\Symp_0(M,\om),\id)
\to H^1(M,\R)$$ can be defined as follows. For an element $\phi
\in \pi_1(\Symp_0(M,\om),\id)$ and for
a class $a \in H_1(M,\Q)$ set 
$$
(F_\om(\phi), a) = ([\om],  \partial_{\phi}a),
$$
where $\partial_{\phi}$ is the homomorphism defined in \S1 and $(\cdot,\cdot)$
is the natural pairing. It is well-known that $\phi$ is represented by a Hamiltonian
loop if and only if $F_\om(\phi) = 0$.  (See Chapter 10 of~[MS2], for example.)
\MS

\NI
{\bf Proof of 2.A:} Observe that the homomorphism $\partial : H_1(M,\Q) \to
H_2(M,\Q)$  associated to the loop $\{\phi_t^1\}$ equals that associated 
to $\{\phi_t^2\}$, since the loops are homotopic in the group of diffeomorphisms.
If we know that $\{\phi_t^1\}$ for example is homotopic to a Hamiltonian
loop then  1.A applied to $\{\phi_t^1\}$ implies that $\partial$ vanishes.
Thus the loop $\{\phi_t^2\}$ has zero flux and hence it is Hamiltonian up to
homotopy. \QED

\medskip

Define {\it the flux subgroup}
$\Gamma_\om \in H^1(M,\R)$ as the image of the flux homomorphism.
The importance of this notion is due to the  fact that
$\Ham(M,\om)$ is $C^1$-closed in $\Symp_0(M,\om)$
if and only if 
$\Gamma_\om$ is a discrete subgroup of $H^1(M,\R)$.
The statement that $\Gamma_\om$ is  discrete is known as {\it the $C^1$-flux
conjecture}. It is proved in various cases  with the use of both ``soft" and
``hard" methods, however it is still unsolved in  full generality: see
Lalonde--McDuff--Polterovich~[LMP1]. The technique of the present paper
allows it to be established in the following new case.

\proclaim Theorem 2.C. The flux conjecture holds
for all closed symplectic manifolds $(M^{2n},\om)$
with  first Betti number equal to 1.

Note that there are plenty of closed  symplectic manifolds
with  first Betti number equal to 1 (see Gompf~[G]), though none of the interesting
new examples are known to have nontrivial $\pi_1(\Symp_0(M);\id)$.
Theorem 2.C follows immediately from the next more general statement.

\proclaim Theorem 2.D. The   rank over $\Z$ of the group $ 
\pi_1(\Symp_0(M))/ \pi_1(\Ham(M))$ (which is identified with $\Ga_\om$
by the Flux homomorphism) is not greater than the first
Betti number of $M$. In particular, it is finitely generated over $\Z$.

\proof{} If the first statement does not hold, there are symplectic loops 
$\phi_1, \ldots, \phi_m$ with $m > \be_1(M)$ 
whose fluxes $\la_i = F_\om(\phi_i)$ are
independent over $\Z$ in $H^1(M, \R)$. Perturb the form $\om$ to a
rational form $\om'$ and then perturb the loops $\phi_i$ 
to $\om'$-symplectic loops $\phi_i'$ with 
fluxes $\la_i'$.  Since the $\la_i'$ are
rational, there is a non-trivial integral 
linear combination of them $\la' = \Si_i n_i
\la'_i$ that vanishes. Therefore by  Theorem 1.A the homomorphism
$\partial_{\phi'}$ associated to the loop 
$\phi' = \Pi_i (\phi'_i)^{n_i}$ vanishes.
Thus $\partial_{\phi} = 0$ for the loop $\phi =  \Pi_i (\phi_i)^{n_i}$
and hence this loop has zero flux. But this
means that $\phi$ is in $\pi_1(\Ham(M))$,  contradicting the
hypothesis.
\QED
\MS

\NI
{ \large \bf \S3. Symplectic fibrations over $S^2$}
\MS

   There is a correspondence between loops in the group of symplectic
diffeomorphisms and {\it symplectic fibrations} over $S^2$ with fiber
$(M, \om)$.   By definition a symplectic fibration is a fibration 
 such that the changes
of trivialisation preserve  a given 
symplectic form $\om$ on the fibers.  In
other words, the structure group of the fibration is $\Symp(M)$.  The
correspondence is  given by assigning 
to each symplectic loop $\phi_{t \in
[0,1]}$ in $\Symp_0(M)$ the fibration $(M, \om) \to P_{\phi} \to S^2$
obtained by gluing a copy of $D_2^+ \times M$ 
with another $D_2^- \times M$
along their boundary in the following way:  
$$
(2\pi t, x) \mapsto  (-2\pi t, \phi_t(x)).
$$
(Here $D_2$ is the closed disc of radius $1$ of the plane.)
In what follows we
always assume that the base $S^2$ is oriented, and with orientation induced
from $D_2^+$. Note that this correspondence can be reversed: given a
symplectic fibration over the oriented $2$-sphere together with an
identification of one fiber with $M$, one can reconstruct the homotopy class
of $\phi$.

An important topological tool for the study of such fibrations 
is the {\it Wang exact sequence}: 
$$
 ...\to H_{j-1}(M,\Z) \stackrel{{\partial_{\phi}}}{\to} H_j(M,\Z)
\stackrel{i} \to H_j(P_\phi,\Z) \stackrel{\cap [M]} {\to} H_{j-2}(M,\Z) \to ...
$$
  This sequence can be easily
derived
from the exact sequence of the pair $(P_\phi, M)$, where  
 $M$ is identified with a fiber of $P_\phi$.  The important point for us is, of
course, that the boundary map $H_{j-1}(M)\to H_j(M)$ is precisely the
homomorphism $\p_\phi$ that interests us.  Thus $\p_\phi$ vanishes
exactly when the inclusion $i$ is injective or, equivalently, when the
restriction map $\cap [M]$ is surjective.

We say that a symplectic fibration is {\it
Hamiltonian} if the corresponding loop of symplectomorphisms is homotopic
to a Hamiltonian loop. The crucial point is  that {\it $P_\phi$ is
Hamiltonian if and only if the cohomology class of the 
symplectic structure on the fiber extends to a cohomology class on the total
space}.  This is most easily seen if one considers the Wang sequence on
cohomology
$$
 ...\to H^{j+1}(P_\phi)\stackrel{restr}\to  H^{j+1}(M)
\stackrel{{\partial_{\phi}^*}}{\to} H^j(M)  \to H^{j+2}(P_\phi) \to ... $$
where $\p_\phi^*$ denotes the dual of $\p_\phi$,
and notes that $\phi$ is Hamiltonian exactly when $\p_\phi^*([\om]) = 0$.

With this language, Theorem 1.A above is equivalent to the
following statement.

\proclaim Theorem 3.A.  Let $\phi$ be a Hamiltonian loop
on a closed symplectic manifold $(M,\om)$.
Then the homomorphism $i:H_*(M,\Q) \to H_*(P_{\phi},\Q)$
is injective. 

\medskip

The proof of this statement is sketched in the next section.
The formulation of the rigidity phenomenon 2.A in the language 
of symplectic fibrations is especially simple.

\proclaim Theorem 3.B. Let
$\phi$ be a Hamiltonian loop
on a closed symplectic manifold $(M,\om)$.
Consider the connected 
component of $P_{\phi}$ in the space of all symplectic fibrations with fiber
$M$ and base $S^2$ (where the symplectic form on $M$ is allowed to vary).
Then the whole connected component is formed of Hamiltonian fibrations.

\medskip

\proclaim Remark 3.C. \rm{Certain special cases of  3.A
and 3.B are already
known. One of them was pointed out to us by Seidel,
namely when the structure of the symplectic 
fibration $p: P_{\phi} \to S^2$ comes from a K\"ahler
structure on 
the total space $P_{\phi}$ such that the projection $p$ 
is
holomorphic. In this situation 3.A and 3.B follow from
a result due to Deligne which states that the
Leray spectral sequence of $P_{\phi}$ degenerates: 
see Chapter~3.5 in 
Griffiths--Harris~[GH].  Another special case is
when $\phi$ is generated by a
circle action.  In this case, one considers the equivariant cohomology
$H_{S^1}^*(M,\Q)$ that is defined to be the usual cohomology
of the homotopy quotient 
$$
M/\!/S^1 = ES^1\times_{S^1} M,
$$
where $\pi:  ES^1\to BS^1 = \C P^\infty$ is the universal $S^1$-bundle:
see Kirwan~[K].  
It is easy to check that 
the bundle $P_\phi\to S^2$ is just the restriction of
the bundle $M/\!/S^1\to \C P^\infty$ to
 $\C P^1$.  Further, one can check that the vanishing of $\p_\phi$ is
equivalent to the 
degeneration of the spectral sequence for the cohomology
of $M/\!/S^1$, a fact that  is proved by Kirwan in~[K] by using
localization formulas.  Thus 3.A gives an alternative proof of this
degeneration.}

\MS

\NI

\smallskip

Since $p:P_{\phi} \to S^2$ is a {\it Hamiltonian} fibration
it carries a natural deformation class  of symplectic forms
given by the weak coupling construction. 
Recall that {\it the coupling class} $u_\phi \in
H^2(P_{\phi},\R)$ is the (unique) class whose top power vanishes, and whose
restriction to a fiber coincides with the cohomology class of the fiberwise
symplectic structure. Let $\tau$ be a positive generator
of $H^2(S^2,\Z)$. The deformation class above consists of
symplectic forms $\Om$ which represent
the cohomology
class of the form $ u_\phi + \ka\, p^*\tau$  ($\ka >>0$)
and extend the fiberwise symplectic structure.   It is always possible to
choose $\Om$ so that it is a product with respect to the given product
structure near the fibers $M_0$ at $0\in D_2^+$ and $M_\infty$ at 
$0\in D_2^-$: see the proof of Lemma 3.E below.

Besides the coupling class $u_\phi$,
the total space $P_{\phi}$ carries another canonical
second cohomology class 
$$
c_\phi = c_1(TP_\phi^{{\rm vert}})  \in H^2(P_{\phi},\R)
$$
that is defined to be the first Chern
class of the vertical tangent bundle. 

\proclaim Remark 3.D. {\rm The existence of this extension $c_\phi$ 
of the first Chern class $c_1(TM)$ provides
a natural explanation of a phenomena that 
was  first observed by McDuff in~[M1]
and rediscovered by Lupton--Oprea~[LO], 
namely that the flux homomorphism
$F_\om: \pi_1(\Symp_0(M,\om)) \to H^1(M,\R)$ vanishes when 
the symplectic class $[\om]$ is a multiple of $c_1$.}

 Both  classes $u_\phi, c_\phi$ behave well under
compositions of loops.  More precisely, consider two elements $\phi, \psi\in
\pi_1(\Ham(M,\om))$ and their composite $\psi*\phi$.  This can be
represented either by the product $\psi_t\circ\phi_t$ or by the concatenation
of loops.   It is not hard to check that the bundle $P_{\psi*\phi}$ can be
realised as the fiber sum $P_\psi\#P_\phi$ obtained as follows.
Let $M_{\phi,\infty}$ denote the fiber at $0\in D_2^-$ in $P_\phi$ and
$M_{\psi,0}$ the fiber at $0\in D_2^+$ in $P_\psi$. Cut out open product
neighborhoods of each of  these fibers and then glue the
complements by an orientation reversing symplectomorphism of the
boundary.   The resulting space may be realised as
$$
D_2^+\times M\;\cup_{\al_{\phi, -1}}\;S^1\times [-1,1]\;\cup_{\al_{\psi,1}} 
\;D_2^-\times M,
 $$
where
$$
{\al_{\phi, -1}}(2\pi t,x) = (2\pi t, -1,\phi_t(x)),\quad 
{\al_{\psi, 1}}(2\pi t, 1,\psi_t(x)) = (2\pi t,x),
$$
and this may clearly 
be identified with  $P_{\psi*\phi}$. Set
$$
V_\phi = D_2^+\times M\;\cup\;S^1\times [-1,1/2), \quad 
V_\psi = \;S^1\times (-1/2,1]\;\cup 
\;D_2^-\times M.
$$

The next lemma follows imediately from  the construction of the coupling form
via symplectic connections: see~[P2] or [MS2].

\proclaim Lemma 3.E. The classes $u_{\psi*\phi}$ and $c_{\psi*\phi}$ are
compatible with the decomposition $
P_{\psi*\phi}= V_\psi\cup V_\phi$ in the sense that their
 restrictions to $V_\psi\cap V_\phi = (-1/2,1/2)\times S^1 \times M$ equal
the pullbacks of $[\om]$ and $c_1(TM)$.

\proclaim Corollary 3.F.  For every $k \in \{1,...,n\}$ the map
$$
\phi \mapsto \int_{P_\phi} (c_\phi)^k (u_\phi)^{n+1-k}
$$
defines a homomorphism $I_k: \pi_1(\Ham(M,\om))\to \R$.

\NI
{\bf  Remark 3.G.} 
When $(M, \om)$ is monotone
\footnote
{
In this note we will say that $(M,\om)$ is monotone if, for some 
positive $\ka\in
\R$,
 $c_1(TM) = \ka [\om]$ on the whole of $H_2(M)$.  If 
this equation holds only 
on the spherical part $H_2^S(M)$ of $H_2$ we will
call $(M,\om)$ spherically monotone.  
}
the homomorphism $I_1$ agrees with the mixed
action--Maslov homomorphism $I$ defined by 
Polterovich in [P1].  However,
although they are both defined in  
the spherically monotone case, they can differ
since  $I$ depends only on the values of $c_\phi$ on spheres, 
while $I_1$ may not.
Indeed 
$$
I_1(\phi) =  c_\phi\left(\PD((u_\phi)^n\right),
$$
and the Poincar\'e dual $\PD((u_\phi)^n)\in H_2(P_\phi)$ need not be
in $H_2^S(M)\otimes \R$.  For example, if $M$ is a nontrivial $S^2$-bundle
over a Riemann surface of genus $> 0$ and $\phi$ is given by an $S^1$-action
that rotates the fibers of $M$ it is not hard to check that
$\PD((u_\phi)^2)$  is not spherical.

\bigskip

\NI
{ \large\bf \S4 Seidel's maps $\Psi_{\phi,\si}$}
\MS

We start with the definition of the quantum cohomology ring of $M$.
In view of our purposes in the next section, we will give
two versions of this definition,  one with rational
and one with real coefficients.
To simplify our formulas we will
denote the first Chern class $c_1(TM)$ of $M$ by $c$.

Let $\La$ be the usual (rational) Novikov ring of the group $\Hh =
H_2^S(M,\Z)/\!\!\sim$  with valuation $\om(.)$ where $B\sim B'$
if $\om(B-B') = c(B-B') = 0$, and   let $\La_R$ be
the analogous (real) Novikov ring based on the group $\Hh_R =
H_2^S(M,\R)/\sim$. Thus the elements of $\La$ have the form $$
\sum_{B\in \Hh} \la_B e^B
$$
where for each $\ka$ there are only finitely many nonzero
$\la_B\in \Q$ with $\om(B) < \ka$,
 and the elements of
$\La_R$ are
$$
\sum_{B\in \Hh_R} \la_B e^B,
$$
where 
$\la_B\in \R$ and there is a similar finiteness condition.
\footnote
{In [Se] Seidel works with a simplified version of the Novikov
ring $\La$ where the coefficients $\la_B$ affecting $e^B,\; B \in \Hh$,
are elements of $\Z/2\Z$. However, his results extend
in a staightforward way to the case of rational coefficients
by taking into account orientations on the moduli spaces
of pseudo-holomorphic curves. Let us emphasize that in our
definition of $\La_R$ not only the coefficients $\la_B$ are real,
but also the exponents $B$ belong to a real vector space $\Hh_R$.}
Set 
$QH_*(M)
= H_*(M)\otimes\La$ and $QH_*(M,\La_R)
= H_*(M)\otimes\La_R$.  Then $QH_*(M)$ 
is $\Z$-graded with $\deg(a\otimes e^B) = \deg(a) - 2c(B)$.  
It is best to
think of
$QH_*(M,\La_R)$ as $\Z/2\Z$-graded with 
$$
QH_{\ev} =  
H_{\ev}(M)\otimes\La_R, \quad QH_{{\rm odd}} =  
H_{{\rm odd}}(M)\otimes\La_R.
$$
 With respect to the quantum intersection product (defined in \S5
below) both versions of quantum homology are graded-commutative rings
with unit $[M]$. Further, the units in $QH_{\ev}(M,\La_R)$ form a group 
$QH_{\ev}(M,\La_R)^\times$ that acts on   $QH_*(M, \La_R)$ by
quantum multiplication. 

Now we describe how Seidel arrives at an
action of the loop $\phi$ on the quantum homology of $M$.
Denote by $\Ll$ the space of
contractible loops in the manifold $M$.  
Fix a constant loop $x_* \in \Ll$, and
define a covering $\Tilde \Ll$ of $\Ll$ with the base point
$x_*$ as follows. Note first that a path between $x_*$ and
a given loop $x$ can be considered as a $2$-disc $u$ in $M$
bounded by $x$. Then the covering $\tilde \Ll$ is defined by saying that
two paths are equivalent if the $2$-sphere $S$ obtained by gluing
the corresponding discs has $\om(S) = c(S) = 0$. 
Thus the covering group
of $\Tilde \Ll$ coincides with the abelian group $\Hh$.

Let $\phi = \{\phi_t\}$ be a loop of Hamiltonian diffeomorphisms.
Because the orbits $\phi_t(x), t\in [0,1],$ of $\phi$ are contractible
(see~[LMP1]), one can define a mapping $T_{\phi}: \Ll \to \Ll$ 
which takes
the loop $\{x(t)\}$ to a new loop $\{\phi_t x(t)\}$. Let
$\Tilde T_{\phi}$
be a lift of $T_{\phi}$ to $\tilde \Ll$. To choose
such a lift one should specify a homotopy class of paths
in $\Ll$ between the constant loop and an orbit of $\{\phi_t\}$.
It is not hard to see that in the language 
of symplectic fibrations this choice of lift
corresponds to a choice of an equivalence class $\si$ of 
sections of $P_{\phi}$,
where two sections are identified if their values under
$c_\phi$ and $u_\phi$ are equal.
Thus the lift can be labelled $\Tilde T_{\phi,\si}$.

Recall now that the Floer homology $HF_*(M)$  can be
considered as the Novikov homology of the action functional
on $\Tilde \Ll$.
Therefore $\Tilde T_{\phi, \si}$ defines a natural automorphism $(\Tilde
T_{\phi,\si})_*$ of $HF_*(M)$. Further, 
if $\Phi: HF_*(M) \to HQ_*(M)$ is the
canonical isomorphism constructed in Piunikhin--Salamon--Schwartz~[PSS],
there is a corresponding 
automorphism $\Psi_{\phi,\si}$ of $QH_*(M)$ given by
$$
\Psi_{\phi,\si}\; = \;\Phi \circ ({\Tilde T_{\phi,\si}})_* \circ \Phi^{-1}.
$$
This gives rise to an action of the group of all pairs $(\phi,\si)$ on $QH_*(M)$.

Seidel shows that when $M$ satisfies a suitable semi-positivity condition the map
$\Psi_{\phi,\si}: QH_*(M)\to QH_*(M)$ is in fact induced by quantum multipication
by an element of $QH_{\ev}(M)^\times$ that is formed from the moduli space of all
$J$-holomorphic sections of $P_\phi$.  
 In our work we in a sense go backwards. 
 We give a  new
definition of the maps $\Psi_{\phi,\si}$ that does not explicitly mention Floer
homology and will prove that they are isomorphisms by a direct gluing
argument. 
We  will see in the next section that our 
map does agree with Seidel's. Further, our version of the definition
 no longer restricts us to using the 
coefficients $\La$ via the covering $\Tilde
\Ll\to \Ll$.  Instead we will consider the extension $\La_R$, which
will allow us to  define an action of the group $\pi_1(\Ham)$ itself.

\MS

 Let $\Om$ be  a
symplectic form  on $P_\phi$ that extends $\om$ and is in the
natural deformation class $u_\phi + \ka\,p^*(\tau)$.  As above, define an
equivalence relation on the set of homology classes of sections
of $P_{\phi}$ 
by identifying two such classes if their values under $c_\phi$ and $u_\phi$
are equal.   
Then, given a loop of Hamiltonian diffeomorphisms $\phi$ on $M$,
and an equivalence class of sections $\si$ of $P_{\phi}$
with $d = 2c_\phi(\si)$, define a $\La$-linear map
$$
\Psi_{\phi,\si}: \;\; QH_*(M) \to QH_{*+d}(M)
$$
as follows.   For $a \in
H_*(M,\Z)$, $\Psi_{\phi,\si}(a)$ is  the class in $QH_{*+d}(M)$ whose
intersection with $b\in H_*(M,\Z)$ is given by: 
$$
  \Psi_{\phi,\si}(a) \cdot_M
b = \sum_{B\in \cal H} n(i(a),i(b); \si + i(B))e^B.
$$
Here $n(v,w;D)$ denotes the Gromov--Witten invariant
which counts isolated $J$-holomorphic stable curves 
in $P_\phi$ of genus $0$
and two marked points that
represent the equivalence class $D$ and whose marked points go
through  given generic
representatives of the classes $v$ and $w$ in $H_*(P_{\phi},\Z)$. 
More precisely, one  defines $n(v, w; D)$ to be
the intersection of the virtual moduli cycle
$$
\ev: \oMm_{0,2}^\nu(P_\phi, J, D) \to P_\phi\times P_\phi,
$$
that consists of all perturbed 
$J$-holomorphic genus $0$ stable maps that lie in
class $D$ and have $2$ marked points, with a generic representative of
the class
$v \otimes w$ 
in $P_\phi\times P_\phi$. This definition is well understood
provided 
$M$  is spherically monotone or has minimal
spherical  Chern number~\footnote
{
The minimal spherical Chern number $N$
is the smallest positive 
integer such that the image of $c = c_1(TM)$ on
$H_2^S(M)$ is contained in $N\Z$.
It equals to $+\infty$ when this image vanishes.
}
 $N \ge n-1$. In the general case, the definition of Gromov-Witten
invariants along these lines
forms the subject of recent works [FO],[LiT],[LT],[R] and [HS].
Further, we have written $i$ for the homomorphism
$H_*(M) \to H_*(P)$ and $\cdot_M$ for the linear extension to $QH_{\ast}(M)$
of the usual  intersection
pairing on $H_*(M,\Q)$.  Thus
$a\cdot_M b = 0$ unless $\dim(a) + \dim (b) = 2n$ in which case it is the
algebraic number of intersection points of the cycles.  Note finally that,
by Gromov compactness, there are for each
given energy level $\ka$  only finitely many homology classes $D$
with $\om(D-\si) \le \ka$ that are represented by $J$-holomorphic curves
in $P_\phi$.  Thus
 $\Psi_{\phi,\si}(a)$ satisfies the finiteness condition for
elements  of $QH_*(M,\La)$.

Since $n(i(a), i(b); D) = 0$ unless $
2c_\phi(D) + \dim (a) + \dim (b) = 2n,
$
we have 
$$
\Psi_{\phi,\si}(a) = \sum a_B\otimes e^B,
$$
where 
$$
\dim (a_B) = \dim(a) + 2c_\phi(D)  = \dim (a) + 2c_\phi(\si) + 2c(B).
$$
 Observe also that
$$
\Psi_{\phi,\si + B} = \Psi_{\phi,\si}\otimes e^{-B}.
$$

When $M$  is spherically monotone or has minimal
spherical  Chern number at least $n-1$
the following two results are proved by Seidel.  The general
case will be established in~[LMP2].  

\proclaim Lemma 4.A. If $\phi$ is the constant loop $*$ and $\si_0$ is the flat
section $pt\times S^2$ of $P_* = M\times S^2$ then $\Psi_{*,\si_0}$ is the
identity map.

 \proclaim Proposition 4.B.  Given  sections $\si$ of $P_\phi$ and  $\si'$ of
$P_\psi$ let $\si'\#\si$ be the union of these sections in the fiber sum
$P_\psi\#P_\phi = P_{\psi*\phi}$.  Then 
$$
\Psi_{\psi,\si'}\circ \Psi_{\phi,\si} = \Psi_{\psi*\phi,
\si'\#\si}.
$$

The main step in the proof of these statements
is to show that when
calculating the Gromov-Witten invariant $n(i(a),i(b);D)$
via the intersection between the virtual moduli cycle
and the class $i(a) \otimes i(b)$
we can assume the following:

\NI
--- the representative of $i(a) \otimes i(b)$ has the form
$\al\times \be$ where $\al, \be$ are cycles 
lying in the fibers of $P_\phi$;
 
\NI
--- the intersection occurs with elements in the top stratum of
$\oMm_{0,2}^\nu(P_\phi, J, D)$ consisting of sections of $P_\phi$.

\NI
Lemma~4.A is then almost immediate, and  
Proposition~4.B can be proved by the
well-known  gluing techniques of [RT] or [MS1].

\proclaim Corollary 4.C.  $\Psi_{\phi,\si}$ is an isomorphism for all loops
$\phi$ and sections $\si$. 

With this in hand, we can establish
Theorem 3.A and hence also 1.A.
\SS

\NI
{\bf Proof of 3.A:}

Gromov-Witten invariants are linear in each variable.
Thus if $i(a)=0 $ for some $a \neq 0$, then $\Psi_{\phi,\si}(a) = 0$,
a contradiction with the fact that $\Psi_{\phi,\si}$ is an 
isomorphism. \QED
\MS

\NI
{ \large\bf \S5  The representation of $\pi_1(\Ham(M))$}
\MS

In this section we prove the 
following mild generalization of the main result in [Se].

\proclaim Theorem 5.A. There exists a
homomorphism 
$$
\Psi:  \pi_1(\Ham(M,\om)) \to QH_{\ev}(M,\La_R)^\times.
$$

Our homomorphism is obtained from Seidel's by a process of averaging, 
and contains much  the same information.\footnote
{
In fact, when $(M,\om)$ is 
spherically monotone, we may take the range of
this homomorphism to be $QH_{\ev}(M,\La)^\times$.   This case of our
theorem was proved in the first version of Seidel's paper.}
 In particular, his
calculations show that it is nontrivial in many cases. 
Our averaging procedure forces us to work with the {\it real}
Novikov ring $\La_R$ which was introduced in the previous section.
Note also that one 
cannot always replace $\La_R$ by $\La$ even when $\om$ 
is integral unless $(M,\om)$ is spherically monotone.

In order to use the maps
$\Psi_{\phi,\si}$ to define a representation of the 
group $\pi_1(\Ham(M,\om))$
we must  make a canonical choice of 
section $\si_\phi$ that (up to  equivalence)
satisfies the composition rule $$
\si_{\psi*\phi} = \si_\phi\#\si_\psi,
$$
where $\si_\phi\#\si_\psi$ denotes the  obvious union of the sections in
$P_{\psi*\phi} = P_\psi\#P_\phi$.  Unfortunately, it is not always possible to
do this if one just considers usual sections.  Further, one has to proceed
slightly differently in the case when 
 the classes $[\om]$ and $c = c_1(TM)$ are linearly
dependent on $H_2^S(M)$.  So let us assume to begin with that these classes
are linearly independent.

We will say that $\si$ is an $\R$-section of
$P_\phi$ if it is a finite sum $\sum\la_i\si_i, \la_i\in \R,$ of sections such that
$\sum \la_i = 1$.  Then, by our assumption on
$[\om]$ and $c$, there
is an $\R$-section  $\si_\phi$ such that  $$
 u_\phi(\si_\phi) =
0,\quad c_\phi(\si_\phi) = 0. 
$$ 
Clearly, the equivalence class of this
$\R$-section is unique.  Also, by Lemma~3.E the needed composition rule
holds.  Further, the definition of the map $\Psi_{\phi,\si_\phi}$ still
makes perfect sense provided that one allows the coefficients $B$ to belong
to $\Hh_R = H_2^S(M,\R)/\!\!\sim$ so that the sum $\si_\phi + B$ can be
integral.  We therefore get  a representation 
$$
\rho:\; \pi_1(\Ham(M,\om))\to \Hom_{\La_R}(QH_*(M,\La_R))
$$
of $\pi_1(\Ham(M))$ in the group of  automorphisms of the $\La_R$-module
$QH_*(M,\La_R)$.

One should think of the $\R$-section $\si_\phi$ as an average of the sections in
$P_\phi$.  The effect of enlarging the Novikov ring to $\La_R$ is 
thus to make enough room
to take this average.  

\MS

Now consider the  case when the classes $[\om]$ and $c$ are linearly
dependent on $H_2^S(M)$.  The difficulty here is that 
the canonical extensions $u_\phi$ and $c_\phi$ need not be dependent
on $H_2^S(P_\phi)$.
 (For example, consider the case when 
 $\phi$ is a  rotation of $M = S^2$.)  Therefore, there
may be no $\R$-section such that
 $
 u_\phi(\si_\phi) =
0,$ $ c_\phi(\si_\phi) = 0. $  However, in this case, the equivalence relation on
$\Hh_R$ is given simply by $[\om]$.  Moreover, if the class $u_\phi$ has the
same value on the two sections $\si, \si'$, so does $c_\phi$.  Hence it suffices to
choose $\si_\phi$ so that $u_\phi(\si_\phi) = 0$.  The value of $c_\phi$ on
$\si_\phi$ is the same for all choices of $\si_\phi$   (though it may not be zero),
and so $\si_\phi$ is still unique up to equivalence.  Thus the previous
arguments go through.
\MS

 In order to complete the proof of  Theorem 5.A
we have  to show  that
the automorphisms in the image of $\rho$  commute with the quantum intersection 
product.  To do this, it is useful to
 describe the homomorphism $\rho$ in the terms used by Seidel.
Recall that the quantum intersection product $a*_Mb$ of two classes $a\in H_i(M,
\Q)$, $b\in H_j(M)$ is defined as follows:
$$
a*_Mb = \sum_{B\in \Hh} (a*_Mb)_B\otimes e^B,
$$ 
where  $(a*_Mb)_B\in H_{i+j- 2n+2c(A)}$ is defined by the requirement that
$$
(a*_Mb)_B\,\cdot_M\, c = n_M(a,b,c;B),
$$
where $n_M(a,b,c;B)$ is the ``number of isolated $J$-holomorphic spheres in class
$B$ that meet $a,b,$ and $c$".  More precisely, $n_M(a,b,c;B)$ is the 
Gromov--Witten invariant that counts the number of (perturbed) $J$-holomorphic
curves in class $B$ that meet the classes $a,b$ and $c$.  This product is extended to
$QH_*(M)$ by linearity over $\La$.  It clearly extends also to $\La_R$. Note here
that when defining $a*_Mb$ we still sum over classes $B\in \Hh$ (and not $B\in
\Hh_R$), since $J$-holomorphic spheres can only represent integral classes.

The next proposition appears in Seidel when $M$ satisfies the semi-positivity
condition described in \S6.  The proof in general follows by looking at what
happens to the Gromov--Witten invariants
$n(i(a), i(b), D)$ when the representatives of $i(a)$ and $i(b)$ are taken to lie in
the same fiber.

 \proclaim  Proposition 5.B.  For all $\phi\in \pi_1(\Ham(M))$,
$$
\rho(\phi)(a) = \Psi_{\phi,\si_\phi}([M])\,*_M\; a.
$$

\NI
{\bf Proof of Theorem 5.A.}

This is now clear.  By Proposition 5.B the homomorphism 
$$
\Psi: \quad \pi_1(\Ham(M)) \to
QH_*(M;\La_R)^\times
$$ 
is given by
$
\phi\mapsto \rho(\phi)([M]).
$
\QED
\MS

\NI
{ \large\bf \S6 What is proved and what will be proved}
\MS

The outline of the proof of Theorem 1.A 
(or of the equivalent Theorem 3.A)
given above
becomes completely rigorous provided $(M,\omega)$
satisfies one of the following assumptions which
were used by Seidel in [Se] for definition of his action:

(i) $(M,\omega)$ is spherically monotone;

(ii) The minimal spherical Chern number $N$ does not exceed
$n-1$.

\NI
The same is true for our results in \S 5.

The reduction of our results in section 2 and 3 
to  1.A is ``soft" and works without any
additional conditions. The only point where
one should be careful is that we need 1.A to hold
{\it simultaneously}
for all symplectic forms under consideration.
In view of this,
2.A is proved provided each of two symplectic forms
satisfy either (i) or (ii), while for 2.B and 3.B one
needs to assume the deformation invariant assumption (ii).

In [LMP2] we will prove 1.A in full generality,
and thus all the results of the present paper
will be confirmed without additional assumptions.

\bigskip

\NI
{\bf Acknowledgment.} We thank Paul Seidel for 
stimulating discussions and for useful
critical remarks on a preliminary version of this paper. We are also
grateful to Dietmar Salamon for helpful discussions.

\bigskip

\NI
{\bf References}
\bigskip

\NI
[BP] M. Bialy and L. Polterovich, Hamiltonian
diffeomorphisms and Lagrangian distributions,
{Geometric and Funct. Analysis}, {\bf 2} (1992), 173-210.
\MS

\NI
[FO]  K. Fukaya and K. Ono, Arnold conjecture and
Gromov--Witten invariants, preprint (1996)
\MS

\NI
[G]  R.~Gompf,
    A new construction for symplectic $4$-manifolds, 
    {\it Annals of Mathematics}, {\bf 142} (1995), 527--595.
\MS

\NI
[GH]
P. Griffiths and J. Harris, {\it Principles of algebraic geometry}.
Wiley, New York (1978). 
\MS

\NI
[HS] H. Hofer and D. Salamon, in preparation.
\MS

\NI
[K]  F.
Kirwan,  {\it Cohomology of quotients in symplectic and algebraic
geometry\/}. Mathematics Notes, {\bf 31}. Princeton University Press
(1984). 
\MS

\NI
[LMP1] F. Lalonde, D. McDuff and L. Polterovich, On the Flux
conjectures, preprint dg-ga/9706015,
 to appear in the Proceedings of the CRM Workshop
on Geometry, Topology and Dynamics, 
Montreal 1995, CRM Special Series pubished by the AMS, 1997.
\MS

\NI
[LMP2] F. Lalonde, D. McDuff and L. Polterovich, 
in preparation.
\MS

\NI
[LiT]  Jun Li and G. Tian, Virtual moduli cycles and Gromov--Witten invariants of
general symplectic manifolds, preprint (1996)
\MS

\NI
[LT] G. Liu and G. Tian, Floer homology and Arnold
conjecture, preprint (1996)
\MS

\NI
[LO]
G. Lupton and J. Oprea,  Cohomologically symplectic spaces, 
Toral actions and the Gottlieb group. Preprint.
\MS

\NI
[M1] 
D. McDuff,   Symplectic diffeomorphisms 
and the flux homomorphism, {\it Inventiones Mathematicae}, {\bf 77},
(1984) 353--66. 
\MS

\NI
[M2] 
D. McDuff,  Recent developments in Symplectic Topology, to appear in {\it
Proceedings of the
2nd European Congress}, Budapest (1996)
\MS

\NI
[MS1]   D. McDuff and D.A. Salamon, {\it $J$-holomorphic curves and
quantum cohomology}, Amer Math Soc Lecture Notes \#6, Amer.
Math. Soc. Providence (1995).
\MS

\NI
[MS2]  D. McDuff and D.A. Salamon, {\it Introduction to
Symplectic Topology},  OUP, Oxford, (1995)
\MS

\NI
[PSS] S. Piunikhin, D. Salamon and M. Schwarz, Symplectic Floer--Donaldson
theory and Quantum Cohomology, {\it Contact and Symplectic Geometry} ed C.
Thomas, Proceedings of the 1994 Newton Institute Conference, CUP, Cambridge (1996)
\MS

\NI
[P1]
L. Polterovich, Hamiltonian Loops and Arnold's 
principle, preprint (1996), to appear in the volume in honor
of V.I. Arnold.
\MS

\NI
[P2]
L. Polterovich, Symplectic aspects of the first eigenvalue, 
preprint dg-ga/9705003.
\MS

\NI
[R]
Y. Ruan, Virtual neighbourhoods and pseudo-holomorphic
curves, Preprint alg-geom/9611021
\MS

\NI
[RT]
Y. Ruan and G. Tian,  A mathematical theory of quantum cohomology.
 Journ Diff Geo {\bf 42} (1995), 259--367.
\MS

\NI
[Se]  P. Seidel, $\pi_1$  of symplectic automorphism groups
and invertibles in quantum cohomology rings, preprint dg-ga/9511011,
(1995) revised (1997), to appear in {\it Geom. and Funct. Anal.}

\end{document}